\documentclass[12pt]{amsart}

\usepackage{amsmath}
\usepackage{amsthm}
\usepackage{amssymb}
\usepackage{amscd}

\newcommand{\gothic}{\mathfrak}

\newcommand{\m}{{\gothic{m}}}
\newcommand{\n}{{\gothic{n}}}
\newcommand{\con}{{\gothic{c}}}
\newcommand{\End}{\operatorname{End{}}}
\newcommand{\Hom}{\operatorname{Hom{}}}
\newcommand{\rank}{\operatorname{rank{}}}
\newcommand{\depth}{\operatorname{depth}}

\newcommand{\syz}{\operatorname{syz{}}}

\renewcommand{\hat}{\widehat}
\renewcommand{\bar}{\overline}

\renewcommand{\phi}{\varphi}

\renewcommand{\to}{{\longrightarrow}}

\newcommand{\charac}{\operatorname{char}}
\newcommand{\e}{\operatorname{e}}

\newtheorem{thm}{Theorem}
\newtheorem{cor}[thm]{Corollary}
\newtheorem{prop}[thm]{Proposition}
\newtheorem{lemma}[thm]{Lemma}

\begin{document}

\title{Hypersurfaces of Bounded Cohen--Macaulay Type}

\author{Graham J. Leuschke}

\address{Department of Mathematics \\
        University of Kansas \\
        Lawrence, KS
        66045}

\email{gleuschke@math.ukans.edu}

\urladdr{http://www.math.ukans.edu/\textasciitilde gleuschke}

\author{Roger Wiegand}

\address{Department of Mathematics and Statistics\\
     University of Nebraska--Lincoln\\
    Lincoln, NE
    68588-0323 }

\email{rwiegand@math.unl.edu}
\urladdr{http://www.math.unl.edu/\textasciitilde rwiegand}

\date{\today}

\thanks{Leuschke's research was supported
by an NSF Postdoctoral Fellowship, and Wiegand's was supported by
grants from the NSA and the NSF}

\bibliographystyle{alpha}

\numberwithin{thm}{section}

\begin{abstract}  Let $R = k[[x_0,\dots,x_d]]/(f)$, where $k$ is a
field and $f$ is a non-zero non-unit of the formal power series
ring $k[[x_0,\dots,x_d]]$.  We investigate the question of which
rings of this form have bounded Cohen--Macaulay type, that is,
have a bound on the multiplicities of the indecomposable maximal
Cohen--Macaulay modules.  As with
 {\it finite} Cohen--Macaulay type, if the characteristic
is different from two, the question reduces to the one-dimensional
case:  The ring $R$ has bounded Cohen--Macaulay type if and only
if $R \cong k[[x_0,\dots,x_d]]/(g+x_2^2+\dots+x_d^2)$, where $g\in
k[[x_0,x_1]]$ and $k[[x_0,x_1]]/(g)$ has bounded Cohen--Macaulay
type.  We determine which rings of the form $k[[x_0,x_1]]/(g)$
have bounded Cohen--Macaulay type.

\vskip 20pt

\noindent This paper is dedicated to our friend and colleague
Wolmer Vasconcelos.

\end{abstract}

\maketitle
\setcounter{section}{-1}
\section{Introduction}

Throughout this paper $(R,\m,k)$ will denote a Cohen--Macaulay local ring
(with maximal ideal $\m$ and residue field $k$).  A maximal Cohen--Macaulay
$R$-module (MCM module for short) is a finitely
generated $R$-module with $\depth(M) = \dim(R)$.
We say that $R$ has {\it bounded} Cohen--Macaulay (CM) type provided
there is a bound on the multiplicities of the indecomposable MCM modules.
One goal of this paper is to examine the distinction between this property
and the formally stronger property of {\it finite} CM type---that there exist,
up to isomorphism, only finitely many indecomposable MCM modules.

We denote the multiplicity of a finitely generated module  $M$ by
$\e(M)$. Following Scheja and Storch \cite{Scheja-Storch}, we say
that $M$ {\it has a rank} provided $K\otimes_RM$ is $K$-free,
where $K$ is the total quotient ring of $R$ (obtained by inverting
all non-zero-divisors).  If $K\otimes_RM \cong K^r$ we say
$\rank(M) = r$.  In this case $\e(M) = r\cdot\e(R)$,
\cite[(4.6.9)]{BH}, so for modules with rank, a bound on
multiplicities is equivalent to a bound on ranks.

The one-dimensional Cohen--Macaulay local rings of finite CM type have
been completely characterized.  To state the characterization,
we let $\hat R$ be the $\m$-adic completion of $R$ and $\bar R$
the integral closure of $R$ in its total quotient ring $K$.

\begin{thm}\label{Drozd-Roiter} Let $(R,\m)$ be a one-dimensional Cohen--Macaulay local ring.
The following are equivalent:
\begin{enumerate}
\item $R$ has finite Cohen--Macaulay type.
\item $R$ has bounded Cohen--Macaulay type and $\hat R$ is reduced.
\item $R$ satisfies the ``Drozd--Ro\u\i ter conditions":
\begin{enumerate}
    \item[({\bf dr1})]$\e(R) \le 3$; and
    \item[({\bf dr2})]$\frac{\m \bar R + R}{R}$ is cyclic as an
$R$-module.\end{enumerate}\end{enumerate}\end{thm} This result was
asserted (in different but equivalent form)  in a 1967 paper
\cite{Drozd-Roiter} by Drozd and Ro\u\i ter, and they sketched a
proof in the ``arithmetic'' case, where $R$ is a localization of a
module-finite ${\mathbb Z}$-algebra. Their proof that (2)
$\implies$ (3) goes through in the general case, but their proof
that (3) $\implies$ (1) is rather obscure even in the arithmetic
case.  In 1978 Green and Reiner \cite{Green-Reiner} gave detailed
matrix reductions proving that (3) $\implies$ (1) in the
arithmetic case. In \cite{Wiegand:1989} R. Wiegand used the
approach in \cite{Green-Reiner} to prove the theorem under the
additional hypothesis that the residue field $R/\m$ is perfect.
Later, in \cite{Wiegand:1994}, he showed that the theorem is true
as long as the residue field does not have characteristic $2$.
Finally, in his 1994 Ph.D. dissertation \cite{Cimen:thesis}, N.~\c
Cimen completed the intricate matrix reductions necessary to prove
the theorem in general.  In \cite{CWW} one can find a streamlined
proof of everything but the matrix reductions, which appear in
\cite{Cimen:paper}. (A word about the implication (1) $\implies$
(2) is in order.  If $R$ has finite CM type, so has the completion
$\hat R$, \cite[Cor. 2]{Wiegand:1994}.  Now by \cite[Prop.
1]{Wiegand:1994} $\hat R$ is reduced.  Also, we mention that
condition ({\bf dr2}) in (3) implies that $\bar R$ is finitely generated
as an $R$-module; therefore ({\bf dr1)} could be replaced by the
condition that $\bar R$ can be generated by three elements as an $R$-module.)

Thus, for one-dimensional analytically unramified local
Cohen--Macaulay rings, finite CM type and bounded CM type are
equivalent.  A statement of this form --- that bounded representation type implies finite representation type --- is often called the ``first Brauer--Thrall conjecture''; see \cite{Ringel:Report} for some history on this and related conjectures.  In particular, the statement for finite-dimensional algebras over a field is a theorem due to Ro\u\i ter \cite{Roiter:BT1}.  The first example showing that that the two concepts
are not equivalent in the context of MCM modules was given
by Dieterich in 1980 \cite{Dieterich:1980}: Let $k$ be a field of
characteristic $2$, let $A = k[[x]]$, and let $G$ be the
two-element group.  Then the group ring $AG$ has bounded CM type.
Note that $AG \cong k[[x,y]]/(y^2)$.  Thus condition (2) of
Theorem~\ref{Drozd-Roiter} fails, and $AG$ has infinite CM type. Later
Buchweitz, Greuel and Schreyer \cite{BGS} noted that
$k[[x,y]]/(y^2)$ has bounded CM type for {\it every} field $k$. In
\S2 we will show, by adapting an argument due to Bass
\cite{Bass:ubiquity}, that every one-dimensional Cohen--Macaulay
ring of multiplicity $2$ has bounded CM type.

In \S1 we show that, for complete equicharacteristic
hypersurfaces, the question of bounded CM type reduces to the case
of plane curve singularities, and in \S2 we examine that case.  We
are able to answer the question completely:  A one-dimensional
complete equicharacteristic hypersurface $R$ has bounded CM type
if and only if either (a) $R$ has finite CM type, (b) $R \cong
k[[x,y]]/(y^2)$ or (c) $R \cong k[[x,y]]/(xy^2)$. The
indecomposable MCM modules over $k[[x,y]]/(xy^2)$ were classified
by Buchweitz, Greuel and Schreyer (\cite{BGS}, see also
Theorem~\ref{Dinfinity} below). When $k = {\mathbb C}$, the field
of complex numbers, the rings in (b) and (c) are the exactly the
rings of countably infinite CM type discussed in \cite{BGS}. They
are the limiting cases A$_\infty$ and D$_\infty$ of the A$_n$ and
D$_n$ singularities $k[[x,y]]/(x^{n+1}+y^2)$ and
$k[[x,y]]/(x^{n-1} + xy^2)$, both of which have finite CM type. We
note that the families of ideals exhibiting uncountable
deformation type in (3.5) of \cite{BGS} do not give rise to
indecomposable modules of large rank; thus there does not seem to
be a way to use the results of \cite{BGS} to demonstrate unbounded
CM type in the cases not covered by (a), (b) and (c).

\section{Complete equicharacteristic hypersurfaces of
dimension two or more}

 Our goal in this section is to prove the following
analog, for bounded CM type, of a beautiful result of Buchweitz,
Greuel, Kn\"orrer and Schreyer (\cite{BGS}, \cite{Knorrer}) on finite CM type:

\begin{thm} Let $k$ be a field, and let $R =
k[[x_0,\dots,x_d]]/(f)$, where $f$ is a non-zero non-unit of the
formal power series ring $k[[x_0,\dots,x_d]], d\ge 2$.  Assume
 that the characteristic of $k$ is different from $2$.
Then $R$ has bounded CM type if and only if $R \cong
k[[x_0,\dots.x_d]]/(g+x_2^2+\dots+x_d^2)$, for some $g\in
k[[x_0,x_1]]$ such that $k[[x_0,x_1]]/(g)$ has bounded CM type.\end{thm}

\subsection*{The Double Branched Cover}

We first set some notation.  In this section, we will set
$S:=k[[x_0, \dots, x_d]]$, a ring of formal power series over a
field $k$, and we denote its maximal ideal by $\n$. We fix a
nonzero element $f \in \n^2$ and set $R=S/(f)$.  We define the
{\it double branched cover} $R^\sharp$ of $R$ by
$R^\sharp=S[[z]]/(f+z^2),$ where $z$ is a new indeterminate over
$S$. Note that there is a natural injection $R \hookrightarrow
R^\sharp$ and a natural surjection $R^\sharp \to R$ defined by
$\bar{z} \mapsto 0$, where $\bar{z}$ is the coset of $z$ in
$R^\sharp$.  There are functors from the category of MCM
$R$-modules to that of MCM $R^\sharp$-modules, and inversely,
defined as follows.  For a MCM $R$-module $M$, set $M^\sharp =
\syz_{R^\sharp}^1(M),$ and for a MCM $R^\sharp$-module $N$, set
$\bar{N} = N/zN.$ We have the following relation on the
compositions of these two functors:

\begin{prop}\label{knorrer}\cite[Proposition 12.4]{Yoshino:book} With
notation as above, assume that $M$ has no non-zero free summand.
Then $\bar{M^\sharp} \cong M \oplus \syz_R^1(M)$.  If $\charac(k)
\neq 2$, then $\bar{N}^\sharp \cong N \oplus
\syz_{R^\sharp}^1(N).$
\end{prop}

(The restriction on the characteristic of $k$ does not appear in
Yoshino's version, which treats only characteristic zero, but the
proof there is easily seen to apply in this context.)

This allows us to show that bounded CM type ascends to and
descends from the double branched cover. We use two slightly more
general lemmas, the first of which is due to Herzog, and the
second of which says that for a Gorenstein ring $A$, a bound on
the multiplicities of indecomposable MCM $A$-modules is equivalent
to a bound on their numbers of generators. We denote by $\nu_R(M)$
the minimal number of generators required for $M$ as an
$R$-module.

\begin{lemma}\cite[Lemma~1.3]{Herzog:1978}\label{herzog} Let $A$ be a
Gorenstein local ring and let $M$ be an indecomposable nonfree MCM $R$-module.
Then $\syz_A^1(M)$ is indecomposable.\end{lemma}

\begin{lemma}\label{snort} Let $A$ be a Gorenstein local ring.
Then there is a bound on the multiplicities of indecomposable
MCM $A$-modules if and only if there is a bound on the number of
generators of same.\end{lemma}

\begin{proof}  Let $M$ be an indecomposable nonfree MCM $A$-module, $n= \nu_A(M)$, and
let $N = \syz_A^1(M)$, so that we have the following exact sequence
$$\CD 0 @>>> N @>>> A^n @>>> M @>>> 0.\endCD$$
By Lemma~\ref{herzog}, $N$ is also indecomposable.  If the multiplicities
of both $M$ and $N$ are bounded above by $B$, then
$n \e(A) = \e(A^n) \leq 2B$, so $n\leq 2B/\e(A)$.
Conversely, if $n$ is bounded by some number $B$,
then $\e(M) \leq \e(A^n) \leq B\e(A)$, so the multiplicity of
$M$ is bounded.\end{proof}

\begin{prop}\label{sharp}Let $R=S/(f)$ be a complete
hypersurface, where $S=k[[x_0, \dots x_d]]$ and $f$ is a nonzero
nonunit of $S$.
\begin{enumerate}
\item If $R^\sharp$ has bounded CM type, then $R$ has bounded CM
type as well.
\item If the characteristic of $k$ is not $2$, then the
converse holds as well.  In fact, if $\nu_R(M) \le B$ for each
indecomposable MCM $R$-module, then $\nu_{R^\sharp}(N) \le 2B$ for
each indecomposable MCM $R^\sharp$-module $N$.
\end{enumerate}\end{prop}

\begin{proof} Assume (1).  First let $M$ be an indecomposable nonfree MCM $R$-module.
Then by Prop.~\ref{knorrer} $\bar{M^\sharp} \cong M \oplus
\syz_R^1(M)$, so $M$ is a direct summand of $\bar{M^\sharp}$.
Decompose $M^\sharp$ into indecomposable MCM $R^\sharp$-modules,
$M^\sharp \cong N_1 \oplus \cdots \oplus N_t$, where each $N_i$
requires at most $B$ generators. Then $\bar{M^\sharp} \cong
\bar{N_1} \oplus \cdots \oplus \bar{N_t}$, and by the
Krull-Schmidt uniqueness theorem $M$ is a direct summand of some
$\bar{N_j}$.  Since $\nu_R(\bar{N_j}) = \nu_{R^\sharp}(N_j)$, the
result follows.

For the converse, let $N$ be an indecomposable nonfree MCM
$R^\sharp$-module. By Prop. \ref{knorrer}, $\bar{N}^\sharp \cong N
\oplus \syz_{R^\sharp}^1(N)$.  Decompose $\bar{N}$ into
indecomposable MCM $R$-modules, $\bar{N} \cong M_1 \oplus \cdots
\oplus M_s$, with $\nu_R(M_j) \leq B$ for each $j$.  Then
$\bar{N}^\sharp \cong M_1^\sharp \oplus \cdots \oplus M_s^\sharp$.
By the Krull--Schmidt theorem again, $N$ is a direct summand of
some $M_j^\sharp$.  It will suffice to show that
$\nu_{R^\sharp}(M_j^\sharp) \le B$ for each $j$.

If $M_j$ is not free, we have $\nu_{R^\sharp}(M_j^\sharp) =
\nu_R(\bar{M_j^\sharp})=\nu_R(M_j)+\nu_R(\syz^1_R(M))$ by
Prop.~\ref{knorrer}.  But since $M_j$ is a MCM $R$-module, all of
its Betti numbers are equal to $\nu_R(M_j)$,
\cite[(7.2.3)]{Yoshino:book}.  Thus
$\nu_R(\bar{M_j^\sharp})=2\nu_R(M_j)\le 2B$.  If, on the other
hand, $M_j = R$, then $M_j^\sharp = zR^\sharp \cong R^\sharp$, and
$\nu_R(\bar{M_j^\sharp})=1$.
\end{proof}

\subsection*{Multiplicity and Reduction to Dimension One}

Our next concern is to show that a hypersurface of bounded
representation type has multiplicity at most two, as long as the
dimension is greater than one.  This is a corollary of the
following result of Kawasaki (\cite{Kawasaki:1996}, due originally
in the graded case to Herzog and Sanders
\cite{Herzog-Sanders:1988}).  An {\it abstract hypersurface} is a
Noetherian local ring $(R, \m)$ such that the $\m$-adic completion
$\widehat{R}$ is isomorphic to $S/(f)$ for some regular local ring
$S$ and nonunit $f$.

\begin{thm}\label{kawasaki}\cite[Theorem 4.1]{Kawasaki:1996}
Let $(R, \m)$ be an abstract hypersurface of dimension $d$.
Assume that the multiplicity $\e(R)$ is greater than 2.
Then for each $n > e$, the maximal Cohen--Macaulay module
$\syz_R^{d+1}(R/\m^n)$ is indecomposable and
$$\nu_R(\syz_R^{d+1}(R/\m^n)) \geq \binom{d+n-1}{d-1}.$$\end{thm}

\begin{cor}\label{mult} Let $R$ be an abstract hypersurface with
$\dim(R) >1$ and $\e(R) > 2$. Then $R$ does not have bounded CM
type.\end{cor}

\begin{prop}\label{gorp}  Let $(R, \m, k)$ be a
Gorenstein local ring of bounded CM type.
Then $R$ is an abstract hypersurface. If $R$ is complete, $d:=
\dim(R) \geq 2$, and $R$ contains a field of characteristic not
equal to $2$, then either $R$ is regular or ${R} \cong A^\sharp$
for some hypersurface $A$, which also has bounded CM type.
\end{prop}

\begin{proof} To show that $R$ is an abstract hypersurface,
it suffices (as in the proof of \cite[Lemma~1.2]{Herzog:1978}) to
show that the Betti numbers $\beta_R^n(k)$ are bounded. Let $M =
\syz_R^d(k)$ and decompose $M$ into nonzero indecomposable MCM
modules, $M = \bigoplus_{i=1}^t M_i$. Assume that $M_1, \dots,
M_s$ are nonfree and $M_{s+1}, \dots, M_t$ are free modules. By
Lemma~\ref{herzog}, each subsequent syzygy module of $k$ has
exactly $s$ indecomposable summands, so if $B$ is a bound for the
number of generators of MCM modules, then $\beta_R^n(k) \leq sB$
for all $n \geq d$.

For the final statement, we assume that $R$ is complete and not
regular. By the first part, we can write ${R}\cong S/(f)$ where
$S=k[[x_0, \ldots, x_d]]$ is a power series ring over a field and
$d \geq 2$. Write $f = \sum_{i=0}^\infty f_i$, where each $f_i$ is
a homogeneous polynomial in $x_0, \dots, x_d$ of degree $i$.
Since, by Corollary~\ref{mult}, $\e(R) = 2$, we have $f_0=f_1=0$
and $f_2 \neq 0$. We may assume after a linear change of variables
that $f_2$ contains a term of the form $cx_d^2$, where $c$ is a
nonzero element of $k$. Now consider $f$ as a power series in one
variable, $x_d$, over $S' := k[[x_0, \dots, x_{d-1}]]$. As such,
the constant term and the coefficient of $x_d$ are in the maximal
ideal of $S'$. The coefficient of $x_d^2$ is of the form $c+g$,
where $g$ is in the maximal ideal of $S'$. Therefore, by
(\cite[Theorem 9.2]{Lang}), $f$ can be written uniquely in the
form
$$f(x_d) = u(x_d^2 + b_1 x_d +b_2),$$
where the $b_i$ are elements of the maximal ideal of $S'$ and $u$ is a unit of $S$.

We may ignore the presence of $u$, as it does not change $R$.
Then, since $\charac(k) \neq 2$, we can complete the square and,
after a linear change of variables, write
$f = x_d^2 + h(x_0, \dots, x_{d-1})$ for some power series $h\in S'$.
By Prop.~\ref{sharp}, $A := S'/(h)$ has bounded CM type. \end{proof}

\section{Dimension one}\label{dimensionone}

 The results of the
previous section reduce our problem to the case of one-dimensional
hypersurface rings. In this section we will deal with this case.
Some of our results go through for more general one-dimensional
Cohen--Macaulay local rings, not just hypersurfaces.  We note that
over a one-dimensional CM local ring the MCM modules are exactly
the non-zero finitely generated torsion-free modules.

\subsection*{Multiplicity two}

We begin with a positive result, which puts the examples of
\cite{Dieterich:1980} and \cite{BGS} mentioned earlier into a
general context. In the analytically unramified case, the result
below is due to Bass \cite{Bass:ubiquity}. In \cite{Rush:1991}
Rush proved the result in the analytically ramified case, but only
for modules with rank. Here we will show how to remove this
restriction.

\begin{thm}\label{multtwo} Let $(R,\m)$ be a one-dimensional Cohen--Macaulay
local ring with $\e(R) = 2$.  Then every MCM $R$-module is
isomorphic to a direct sum of ideals of $R$.  In particular, every
indecomposable MCM $R$-module has multiplicity at most $2$ and is
generated by at most $2$ elements.
\end{thm}

\begin{proof} We note that every ideal of $R$ is generated by two
elements,
\cite[Chap.~3,~Theorem 1.1]{Sally}.  If the integral closure
$\bar R$ of $R$ in the total quotient ring $K$ of $R$ is
 finitely generated over $R$, the
theorem follows from \cite[(7.1), (7.3)]{Bass:ubiquity}.  Therefore we
assume from now on that $\bar R$ is {\em not} a finitely generated
$R$-module.

Suppose $S$ is an arbitrary module-finite
extension of $R$ contained in $K$.  By
\cite[Theorem 3.6]{Sally-Vasconcelos:1974}, $\bar R$ is quasi-local,
and it follows that $S$ is local.  Moreover, each ideal of $S$ is
isomorphic to
an ideal of $R$ and is therefore generated by two elements (as an $R$-
or $S$-module). By \cite[(6.4)]{Bass:ubiquity} $S$ is Gorenstein.  
In particular, $R$ itself is Gorenstein.

To complete the proof, it will suffice to show that every MCM $R$-module
$M$ has a direct summand isomorphic to a non-zero ideal of $R$.  The first part
of the argument here is due to Bass \cite[(7.2)]{Bass:ubiquity}.
Suppose first that $M$ is faithful. Let $S = \{\alpha \in K \mid \alpha M \subseteq M\}$.
 Since $M$ is faithful, $S$ is a subring of $\Hom_R(M,M)$ and
 therefore is a module-finite extension
of $R$. Of course $M$ is a MCM $S$-module.  If there is a surjection $M \to S$,
then $M$ has a direct summand isomorphic to $S$, which, in turn, is isomorphic to
an ideal of $R$.
Therefore we suppose to the contrary that $M^* = \Hom_S(M,\n)$, where
$(\_)^*$
denotes the $S$-dual and $\n$ is the maximal ideal of $S$. Now $M^*$ is
a module
over $E:= \Hom_S(\n,\n)$ and therefore so is $M^{**}$.  But since $S$ is Gorenstein,
$M$ is reflexive \cite[(6.2)]{Bass:ubiquity}.  Therefore $M$ is actually an
$E$-module.
Since $S\subseteq E \subseteq K$ we must have $S = E$ (by the definition
of
$S$).  It follows easily that $\n$ is a principal ideal, that is, $S$ is
a
discrete valuation ring. But then $M$ has $S$ as a free summand,
a contradiction.

If $M$ is not faithful, let $I = (0:M)$.  We claim that $M$ has a direct
summand
isomorphic to an ideal of $R/I$.  Since $R/I$ embeds in a direct sum of
copies
of $M$, $R/I$ has depth $1$ and therefore is a one-dimensional
Cohen--Macaulay
ring.  Also, $\e(R/I) \le 2$ since ideals of $R/I$ are two-generated.
Our claim
now follows from the argument above, applied to the $R/I$-module $M$.
To complete the proof, we show that $R/I$ is isomorphic to an ideal of
$R$.
Taking duals over $R$, we note that $(R/I)^* \cong (0:I)$, which, since
$I\ne
0$, is an ideal of height $0$ in $R$.  Therefore $R/(0:I)$ has positive
multiplicity, and by \cite[Chap. 3, Theorem 1.1]{Sally}  $(0:I)$ is a
principal
ideal, that is, $(R/I)^*$ is cyclic.  Choosing a surjection
$R^* \twoheadrightarrow (R/I)^*$ and dualizing again, we have (since
$R/I$ is MCM
and $R$ is Gorenstein) $R/I \hookrightarrow R$ as desired.\end{proof}

\subsection*{Multiplicity at least four}\label{multfoursection}
Next we will show, in (\ref{multfour}), that if $\e(R) \ge 4$ then
$R$ has indecomposable MCM modules with arbitrarily large
(constant) rank. In \cite{Drozd-Roiter} Drozd and Ro\u\i ter
developed a machine for building big indecomposable modules over
certain one-dimensional rings.  Their approach was refined and
generalized in \cite{Green-Reiner} and \cite{Wiegand:1989}.  The
results in \cite{Wiegand:1989} apply to one-dimensional
analytically unramified local rings and use the conductor square
associated to the inclusion $R \hookrightarrow \bar R$, where
$\bar R$ is the integral closure of $R$ in its total quotient
ring.  Here we observe that most of the theory goes through in the
current setting.

\subsection*{Notation and Assumptions}\label{notation}

As always, we assume that $(R,\m)$ is a Cohen--Macaulay local ring
of dimension one with total quotient ring $K$.
 We let $S$ be
a finite birational extension of $R$; that is, $S$ is a subring of $K$
containing
$R$, and $S$ is finitely generated as an $R$-module.
Let $\con$ be the conductor of
$R$ in $S$, that is, the largest ideal of $S$ that is contained in $R$.

We form the conductor square:
\begin{equation}\label{conductorsquare}
\CD R   @>\subset>> S\\
    @VVV    @VV\pi V\\
    \frac{R}{\con} @>\subset>> \frac{S}{\con}
    \endCD\end{equation}

The bottom line of the square $\frac{R}{\con} \hookrightarrow
\frac{S}{\con}$ is an {\it Artinian pair} in the terminology of
\cite{Wiegand:1989}.  (By definition, an Artinian pair is a module-finite
extension of commutative Artinian rings.)  A {\it module} over the
Artinian pair $A\hookrightarrow B$ is a pair $(V,W)$ where $W$ is
a finitely generated projective $B$-module, $V$ is an
$A$-submodule of $W$, and $BV = W$. Morphisms and direct sums are
defined in the obvious way.  We say that the $(A\hookrightarrow
B)$-module $(V,W)$ has {\it constant rank} $r$ provided $W$ is a
free $R$-module of rank $r$.

\begin{lemma}\label{lifting} Suppose the Artinian pair
$\frac{R}{\con}\hookrightarrow \frac{S}{\con}$ in
(\ref{conductorsquare}) has an indecomposable module $(V,W)$ of
constant rank  $r$.  Then there is an indecomposable MCM
$R$-module of constant rank $r$.\end{lemma}

\begin{proof} Let $P$ be a free $S$-module of constant rank $r$
mapping onto $W$ by change of rings.  Define $M$ by the pullback
diagram
\begin{equation}\label{Mpullback}
\CD
    M @>>> P\\
    @VVV   @VV\pi V\\
    V @>\subset>> W
    \endCD
\end{equation}

As in \cite{Wiegand:1989} one checks that (\ref{Mpullback}) is isomorphic to the
pullback diagram

\begin{equation}
\CD
    M @>\subset>>SM\\
    @VVV @VVV\\
    \frac{M}{\con M} @>\subset>> \frac{SM}{\con M}
    \endCD
\end{equation}
where $SM$ is the $S$-submodule of $K\otimes_RM$ generated by the
image of $M$, and the vertical arrows are the natural
homomorphisms.  It follows easily that any non-trivial
decomposition of $M$ would induce a decomposition (non-trivial by
Nakayama's lemma) of $(V,W)$.\end{proof}

We note that in the analytically unramified case the Drozd-Ro\u\i ter
conditions
of Theorem~\ref{Drozd-Roiter} can be translated into conditions on the
the bottom line of the pullback diagram for $R\to \bar R$.  The failure
of these
conditions is exactly what we need to build big indecomposables.

\begin{thm}\label{indecs-artin-pair} Let $A\hookrightarrow B$ be an Artinian pair,
with $(A,\m,k)$ local.  If either
\begin{enumerate}
\item $\nu_A(B) \ge 4$, or
\item $\nu_A(B) = 3$ and $\frac{\m B}{\m}$ is not cyclic as an
$A$-module,
then\end{enumerate}
$A\hookrightarrow B$ has, for each $n$, an
indecomposable module of constant rank $n$.\end{thm}

\begin{proof} Although this result is not stated in the literature, it is
proved in complete detail in \S2 of \cite{Wiegand:1989}. The
context is a bit different there, however, so a brief review of
the proof is in order.  In case (1), it is enough, by
\cite[(2.4)]{Wiegand:1989}, to show that the Artinian pair
$k\hookrightarrow B/\m B$ has big indecomposables of large
constant rank.  The general construction in
\cite[(2.5)]{Wiegand:1989} (with $A := B/\m B$) yields, by
\cite[(2.6)]{Wiegand:1989} and the discussion after its proof, the
desired indecomposables except in the case where $k$ is the
$2$-element field and $B$ has at least $4$ local components. In
this case one can appeal to Dade's theorem \cite{Dade:1963}.

In case (2), we note that $\frac{\m B}{\m}\cong
\frac{A+\m B}{A}$ (since $\m B\cap A = \m$). We put $C = A + \m B$ and
$D = C/\m C$.  By \cite[(2.4)]{Wiegand:1989} it is enough to show that the
Artinian
pair $k\hookrightarrow D$ has big indecomposables of constant rank.  As
shown
in the proof of \cite[(2.4)]{Wiegand:1989}, $D \cong k[x,y]/(x^2,xy,y^2)$, and the
existence of big indecomposables follows from case (i) of \cite[(2.6)]{Wiegand:1989}.
(The
hypothesis that $B$ is a principal ideal ring in \cite[(2.3)]{Wiegand:1989} is
irrelevant
here, as it is used only in the case $\nu_A(B) = 2$.)\end{proof}

\begin{cor}\label{fourgens} With $R$ and $S$ as above, suppose
$\nu_R(S) \ge
4$.  Then, for each $n \ge 0$, there is an indecomposable MCM $R$-module
of
constant $n$.\end{cor}

\begin{proof} Apply Theorem~\ref{indecs-artin-pair} and Lemma~\ref{lifting}.\end{proof}

Finally we are ready to prove our first general result on unbounded CM
type.

\begin{thm}\label{multfour} Let $(R,\m)$ be a one-dimensional Cohen--Macaulay
local
ring with $\e(R) \ge 4$.  Then $R$ has, for each $n$, an indecomposable
MCM
module of constant rank $n$.  In particular, $R$ has unbounded CM type.\end{thm}

By Corollary~\ref{fourgens} it will suffice to prove the following:

\begin{lemma} Let $(R,\m,k)$ be a one-dimensional Cohen--Macaulay
local
ring with $\e(R) = e$.  Then $R$ has a finite birational extension $S$
with
$\nu_R(S) = e$.\end{lemma}

\begin{proof} Let $S_n = \End_R(\m^n) \subseteq K$, and put $S =
\bigcup_nS_n$. To see that this works, we can harmlessly assume
$k$ is infinite.  Let $Rf\subset \m$ be a principal reduction of
$\m$.  Choose $n$ so large that
\begin{enumerate}
\item[(a)]$\m^{i+1} = f\m^i$ for $i \ge n$ and
\item[(b)]$\nu_R(\m^i) = e$ for $i \ge n$.
\end{enumerate}
Since $f$ is a non-zero-divisor (as $R$ is Cohen--Macaulay), it
follows from (a) that $S = S_n$. We claim that $Sf^n = \m^n$.  We
have $Sf^n = S_nf^n\subseteq \m^n$.  For the reverse inclusion,
let $\alpha \in \m^n$. Then $\frac{\alpha}{f^n} \m^n\subseteq
\frac{1}{f^n}\m^{2n} = \frac{1}{f^n}(f^n\m^n) = m^n$.  Therefore
$\frac{\alpha}{f^n} \in S$, and the claim follows.  Now $S$ is
isomorphic to $\m^n$ as an $R$-module, and $\nu(S) = e$ by
(b).\end{proof}

\subsection*{Multiplicity three}

At this point we know
 that a one-dimensional Cohen--Macaulay local ring $R$ has
bounded CM type if $\e(R) \le 2$ and unbounded CM type if $\e(R)
\ge 4$.  Now we address the troublesome case of multiplicity
three for complete equicharacteristic hypersurfaces.

Let $R = k[[x,y]]/(f)$, where $k$ is a field and $f\in (x,y)^3 -
(x,y)^4$.  If $R$ is reduced, we know by (\ref{Drozd-Roiter}) that
$R$ has bounded CM type if and only if $R$ has finite CM type,
that is, if and only if  $R$ satisfies the condition ({\bf dr2}):
$\frac{\m \bar R + R}{R}$ is cyclic as an $R$-module.  If the
characteristic is different from $2,3,5$ there are simple normal
forms \cite{GreuelKnorrer} for $f$, classified by the Dynkin
diagrams D$_n$, E$_6$, E$_7$, E$_8$.  (Of course the A$_n$
singularities, of multiplicity two, have finite CM type too.)
Normal forms have in fact been worked out in all characteristics
\cite{Kiyek-Steinke}, \cite{GreuelKroning}, but the classification
is complicated, particularly in characteristic $2$.  Here we focus
on the case where $R$ is not reduced.

\begin{thm}\label{multthree} Let $R = k[[x,y]]/(f)$, where $k$ is a
field and $f$ is a non-zero non-unit of the formal power series
ring $k[[x,y]]$.  Assume
\begin{enumerate}
\item $\e(R) = 3$.
\item $R$ is not reduced.
\item $R \not\cong k[[x,y]]/(xy^2)$.
\end{enumerate}
For each positive integer $n$, $R$ has an indecomposable MCM module of
constant rank $n$.\end{thm}

The ring $k[[x,y]]/(xy^2)$ does indeed have bounded CM type;
see the discussion following the proof and Theorem~\ref{Dinfinity}.

\begin{proof} We know $f$ has order $3$ and that its
factorization into irreducibles has a repeated factor.
Thus, up to a unit, we have either $f = g^3$ or $f = g^2h$,
where $g$ and $h$ are irreducible elements of $k[[x,y]]$ of
order $1$, and, in the second case, $g$ and $h$ are
relatively prime.  After a change of variables \cite[Cor. 2, p. 137]{Zariski-Samuel}
we may assume that $g = y$.

In the second case, if the leading form of $h$ is not a constant multiple of $y$,
then by \cite[Cor. 2, p. 137]{Zariski-Samuel} we may assume that $h = x$.
This is the case we have ruled out in (3).

Suppose now that the leading form of $h$ is a constant
multiple of $y$.  By a corollary
\cite[Cor. 1, p.145]{Zariski-Samuel} of the Weierstrass Preparation Theorem, there
exist a unit $u$ and a non-unit power series $q \in k[[x]]$ such that
$h = u(y + q)$.  Moreover, $q\in x^2k[[x]]$ (since the leading
form of $h$ is a constant multiple of $y$).
In summary, there are two cases to consider:

\begin{enumerate}
\item{}Case 1: $f = y^3$.
\item{}Case 2: $f = y^2(y +  q)$, $0 \ne q \in
x^2k[[x]]$.
\end{enumerate}

Let $\m$ be the maximal ideal of $R$.  We will show that $R$ has a
finite birational extension $S$ such that $\nu_R(S) = 3$ and
$\frac{\m S}{\m}$ is not cyclic as an $R$-module.  An application of
Lemma~\ref{lifting}
and Theorem~\ref{indecs-artin-pair} will then complete the proof.

In Case (1) we put $S = R[\frac{y}{x^2}] = R + R\frac{y}{x^2} + R
\frac{y^2}{x^4}$.  Clearly $\nu_R(S) = 3$.  It will suffice to show that
$\frac{\m S}{\m^2 S + \m}$ is two-dimensional over $R/\m$.  We
have
$$\m S = \m + R\frac{y}{x} + R\frac{y^2}{x^3} \ \ \ \ \text{and} \ \
\ \
\m ^2 S + \m = \m + R\frac{y^2}{x^2}.$$
We must show (a) $\frac{y}{x} \notin \m + R\frac{y^2}{x^2}$ and
(b) $\frac{y^2}{x^3} \notin R\frac{y}{x} + \m + R\frac{y^2}{x^2}$.
If (a) fails, we multiply by $x^2$ and lift to $k[[x,y]]$, getting
$xy \in (x^3,x^2y,y^2,y^3) \subseteq (x^2,y^2)$, contradiction.
If (b) fails, we multiply by $x^3$ and lift to $k[[x,y]]$,
getting $y^2 \in (x^2y,x^4,x^3y, xy^2,y^3) \subseteq (x,y)^3$,
contradiction. This completes the proof in Case (1).

Assume now that we are in Case (2).
Let $\gothic{A}_r$ denote the set of rings $k[[x,y]]/y^2(y +  q)$
with $q$ an element of order $r$ in $k[[x]]$.  Then  $R \in \gothic
A_r$ for some $r \ge 2$.  (The rings in $\gothic A_1$ are isomorphic
to $k[[x,y]]/(xy^2)$, which we have ruled out.)
Put $z = \frac{y}{x}$, and note that $z^3 + \frac{q}{x}z^2 = 0$.
Therefore $R[z]$ is a finite birational extension of $R$ and
is isomorphic to a ring in $\gothic A_{r-1}$.  If $M$ is any
indecomposable MCM $R[z]$-module, then $M$ is also a MCM $R$-module.
Moreover, any $R$-endomorphism of $M$ is also $R[z]$-linear (as $M$
is torsion-free and $R \to R[z]$ is birational).
It follows that $M$ is indecomposable as an
$R$-module.  Therefore the conclusion of the theorem passes from
$R[z]$ to $R$, and we may assume that $R \in \gothic A_2$.

Thus $R = k[[x,y]]/y^2(y+q)$, where $q$ is an element of order $2$
in $k[[x]]$.  Put $u:= \frac{y}{x^2}$,  $v := \frac{y^2+qy}{x^5}$,
and $S := R[u,v]$.  The relations $u^2 = xv - \frac{q}{x^2}u, uv =
v^2 = 0$ show that $S = R + Ru + Rv$, a finite birational
extension of $R$.  One checks easily that
$$\m S = Rx +
R\frac{y}{x} + R \frac{y^2+qy}{x^4} \ \ \ \ \text{and}\ \ \ \
\m^2 S + \m = Rx + Ry +
R\frac{y^2+qy}{x^3}.
$$
To see that $\nu_R(S) = 3$, it suffices to show that $u \notin R +
\m S$ and $v \notin R + Ru + \m S$.  If $u \in R + \m S$
($ = R + R\frac{y}{x} + R\frac{y^2+qy}{x^4}$), we would
have (after multiplying by $x^4$ and lifting to $k[[x,y]]$)
\begin{equation}\label{whirt}
x^2y \in (x^4,x^3y, y^2+qy, y^2(y+q)).
\end{equation}
In an equation demonstrating this inclusion,
the coefficient of $x^4$ must be divisible by $y$. Cancelling
$y$ from such an equation and combining terms, we get $x^2 =
Ax^3 + B(y+q)$, with $A,B\in k[[x,y]]$.  Writing $q =
Ux^2$ (where $U$ is a unit of $k[[x]]$), we have
$(1-BU)x^2 = Ax^3 + By$.  Since $x^2\notin (x^3,y)$, $B$
must be a unit of $k[[x,y]]$.  But then $y \in (x^2)$,
contradiction.  Suppose now that $v\in R + Ru + \m S $.  Clearing
denominators and lifting to $k[[x,y]]$, we get
\begin{equation}\label{glunk}
y^2+qy \in
(x^5,x^3y,x(y^2+qy), y(y^2+qy)),
\end{equation}
and it follows that $y^2+qy \in
(x^5,x^3y)$.  Proceeding as before, we get an equation $y + q =
Ax^3$. But then $(x^2) = (q) \subseteq (y,x^3)$, contradiction.

Finally, we show that $\frac{\m S}{\m^2 S + \m}$ is two-dimensional over $R/\m$,
that is, (a) $\frac{y}{x} \notin \m^2S + \m$ and (b) $\frac{y^2+qy}{x^4}
\notin R\frac{y}{x} + \m^2S+\m$.  If (a) were false, we could multiply by $x^3$
and get equation (\ref{whirt}), which we have already seen to be impossible.
Suppose now that (b) fails.  Our worn-out argument yields (\ref{glunk}) again,
and the proof is complete.\end{proof}

The argument in the proof of Theorem~\ref{multthree} does not apply
to the ring $R:= k[[x,y]]/(xy^2-y^3)$.  Adjoining the idempotent $\frac{y^2}{x^2}$
to $R$, one obtains a ring isomorphic to $k[[x]]\times
k[[x,y]]/(y^2)$, whose integral closure is $k[[x]] \times
\bigcup_{n=1}^{\infty}R[\frac{y}{x^n}]$.  From this information
one can easily check that $\frac{\m S}{\m}$ is a cyclic $R$-module
for every finite birational extension $S$ of $(R,\m)$, so we cannot apply
Theorem~\ref{indecs-artin-pair}.  We appeal instead to the following result of
Buchweitz, Greuel and Schreyer.

\begin{thm}\label{Dinfinity}\cite[Proposition 4.2]{BGS} Let $P$ be a two-dimensional
regular local ring with maximal ideal $\m$ and let $x, y$ be a
generating set for $\m$. Set $R:= P/(xy^2)$.  Then every
indecomposable MCM $R$-module $M$ has a presentation
$$\CD
0 @>>> P^n @>\phi>> P^n @>>>M @>>>0
\endCD$$
with $n=1$ or $n=2$ and $\phi$ one of the following matrices
\begin{equation*}\begin{split}
n=1\ \ :\ \ & (y), (x), (y^2), (xy), (xy^2) \\
n=2\ \ :\ \ & \left(\begin{matrix} y & x^k\\ 0 & -y \end{matrix}\right),
       \left(\begin{matrix} xy & x^{k+1}\\ 0 & -xy \end{matrix}\right), \\
      & \left(\begin{matrix} xy & x^k\\ 0 & -y \end{matrix}\right),
       \left(\begin{matrix} y & x^{k+1}\\ 0 & -xy \end{matrix}\right)
\end{split}\end{equation*}
where $k = 1, 2, 3....$\end{thm}

\begin{cor} Let $P$ and $R$ be as above.  Then every indecomposable MCM $R$-module
is generated by at most two elements.\end{cor}

\section{Summary}

Let us summarize the results of the previous two sections:

\begin{thm}Let $k$ be any field, and let $R=k[[x_0,\dots,x_d]]/(f)$,
where $d\ge 1$ and $f$ is a non-zero non-unit of the power series
ring $k[[x_0,\dots,x_d]]$.
    \begin{enumerate}
\item{} Suppose $d = 1$ and $R$ is reduced.  Then $R$ has
bounded CM type if and only if $R$ has finite CM type.
\item{} Suppose $d = 1$ and $R$ is not reduced.  Then $R$ has infinite CM type.
$R$ has bounded
CM type if and only if either $R\cong k[[x,y]]/(y^2)$ or $R\cong
k[[x,y]]/(xy^2)$. In more detail:
\begin{enumerate}
    \item{} If $R\cong k[[x,y]]/(y^2)$ or $R\cong
k[[x,y]]/(xy^2)$, then every MCM $R$-module is generated by at
most $2$ elements.
    \item{} If $\e(R) \ge 3$ and $R\not\cong
k[[x,y]]/(xy^2)$, then $R$ has, for each $n\ge 1$, an
indecomposable MCM module of constant rank $n$.
\end{enumerate}
\item{} Suppose $d \ge 2$ and the characteristic of $k$ is
different from $2$.  Then $R$ has bounded CM type if and only if
$R \cong k[[x_0,\dots,x_d]]/(g + x_2^2+\dots+x_d^2)$ for some
$g\in k[[x_0,x_1]]$ for which $k[[x_0,x_1]]/(g)$ is a
one-dimensional ring of bounded CM type.
\end{enumerate}
\end{thm}
\begin{proof}  Items (1) and (3) come from Theorem~0.1 and Prop.~1.5
respectively.  In view of Theorems~2.1, 2.5 and \ref{multthree} we know that, in
the context of item (2), $R$ has bounded CM type if and only if
either $\e(R) = 2$ or $R\cong k[[x,y]]/(xy^2)$.  But if $\e(R) =
2$ then (by an argument like that at the beginning of the proof of
Theorem~\ref{multthree}) one sees easily that $R\cong k[[x,y]]/(y^2)$.
\end{proof}

Using this theorem and the analogous result for finite CM type
\cite{BGS,Knorrer}, one can obtain higher dimensional examples of
rings that have infinite bounded CM type.  For example, take $g =
x_1^2$ in (3) of Theorem 3.1.

We recall \cite{WW:1994} that for one-dimensional local CM rings
of finite CM type there is a universal bound on the multiplicities
of the indecomposable MCM modules.  In fact, if $(R,\m,k)$ is any
one-dimensional local CM ring of finite CM type, then every
indecomposable MCM $R$-module can be embedded in $R^4$ (in $R^3$
if $k$ is algebraically closed).  Since $\e(R) \le 3$, one obtains
a bound of $12$ on the multiplicities of the indecomposable MCM
$R$-modules.  The proof of Lemma~\ref{snort} then gives a crude
bound of $24$ on the number of generators required for the
indecomposable MCM modules.  In fact, the sharp bound on the
number of generators is probably about $8$ and could be determined
by a careful analysis of the work in \cite{WW:1994}.  It is
interesting to observe that one can use Prop.~\ref{sharp} to get such
universal bounds for higher-dimensional hypersurfaces.  Here is a
special case where the sharp bound in dimension one has been
worked out:

\begin{thm} Let $k$ be an algebraically closed field of
characteristic different from $2,3,5$, and let $R\cong
k[[x_0,\dots,x_d]]/(f)$, where $f$ is a nonzero nonunit in the
power series ring.  If $R$ has finite CM type, then every
indecomposable MCM $R$-module can be generated by $6\cdot2^{d-1}$
elements.
\end{thm}
\begin{proof} If $d=1$, one can see from the computations in
Chapter 9 of \cite{Yoshino:book} that every indecomposable MCM
$R$-module is generated by at most $6$ elements.  For $d > 1$ one
uses the main theorem of \cite{BGS,Knorrer}) (the analog of
Prop.~\ref{gorp} for finite CM type) to deduce that $R$ is obtained from
a plane curve singularity of finite CM type by iterating the
``sharp'' operation $d-1$ times.  Then (2) of Prop.~\ref{sharp}
provides the desired bound.
\end{proof}

We have a corresponding result for hypersurfaces of bounded but
infinite CM type. The proof is the same as that of Theorem 3.2. It
is curious that the bound is better than in the case of finite
type. The reason is that by item (2) of Theorem 3.1 the
indecomposable MCM modules in dimension one are generated by two
elements.

\begin{thm}\label{terra} Let $k$ be a field of characteristic not equal to $2$,
and let $R\cong k[[x_0,\dots,x_d]]/(f)$, where $f$ is a nonzero nonunit in the
power series ring.  If $R$ has bounded CM type but not finite CM type, then
every indecomposable MCM $R$-module can be generated by $2^{d}$
elements.\end{thm}


\end{document}